\RequirePackage{fix-cm}
\documentclass[smallextended]{svjour3}       % onecolumn (second format)
\smartqed  % flush right qed marks, e.g. at end of proof
\usepackage{graphicx}

\usepackage{color} 
\def\RR{{I~\hspace{-1.45ex}R} }
\def\NN{{I~\hspace{-1.45ex}N} }

\begin{document}

\title{On the Optimal Control of a Class of Non-Newtonian Fluids%\thanks{Grants or other notes
%about the article that should go on the front page should be
%placed here. General acknowledgments should be placed at the end of the article.}
}
%\subtitle{Do you have a subtitle?\\ If so, write it here}

%\titlerunning{Short form of title}        % if too long for running head

\author{Telma Guerra   \and
        Jorge Tiago \and 
        Ad\'{e}lia Sequeira }

%\authorrunning{Short form of author list} % if too long for running head

\institute{Telma Guerra \at
              Escola Superior de Tecnologia do Barreiro, IPS, CMA, FCT-UNL. \\
                 \email{telma.guerra@estbarreiro.ips.pt}           %  \\
%             \emph{Present address:} of F. Author  %  if needed
           \and
           Jorge Tiago \at
             Departamento de Matem\'{a}tica, CEMAT, Instituto Superior T\'{e}cnico, Universidade de Lisboa\\
              \email{Jorge.tiago{\char'100}dmat.ist.utl.pt}
\and
           Ad\'{e}lia Sequeira \at
              Departamento de Matem\'{a}tica, CEMAT, Instituto Superior T\'{e}cnico, Universidade de Lisboa\\
              \email{adelia.sequeira{\char'100}dmat.ist.utl.pt}              
}

%\dedication{dedication}

\date{Received: date / Accepted: date}
% The correct dates will be entered by the editor

\maketitle

\begin{abstract}
We consider optimal control problems of systems governed by stationary, incompressible generalized Navier-Stokes equations with shear dependent viscosity in a two-dimensional or three-dimensional domain. We study a general class of viscosity functions including shear-thinning and shear-thickening behavior. We prove an existence result for such class of optimal control problems.
% \PACS{PACS code1 \and PACS code2 \and more}
% \subclass{MSC code1 \and MSC code2 \and more}
\end{abstract}

\section{Introduction}
\label{intro}
This paper is devoted to the proof of the existence of  steady solutions for a distributed optimal control problem of a viscous and incompressible fluid. The control and state variables are constrained to satisfy a generalized Navier-Stokes system of equations with shear dependent viscosity which switches from shear-thinning to shear-thickening behavior. More specifically, we deal with the following generalization of the stationary Navier-Stokes system 

\begin{equation}\label{equation_etat}
\left\{ \begin{array}{ll}-div\, \left(S(D\mathbf{\mathbf{y}})\right)+
	\mathbf{\mathbf{y}}\cdot \nabla \mathbf{\mathbf{y}} +\nabla p= \mathbf{u}& \qquad\mbox{in} \ \Omega\vspace{2mm} \\
             div\,\mathbf{\mathbf{y}}=0&  \qquad \mbox{in} \ \Omega\vspace{2mm}\\
             \mathbf{\mathbf{y}}=0& \qquad \mbox{on}
	 \ \partial\Omega\ ,\end{array}\right.\end{equation}
	 
\noindent where $S$ is the extra stress tensor given by

$$S(\eta)=(1+|\eta|)^{\alpha(x)-2}\eta$$

\noindent and $\alpha(x)$ is a positive bounded continuous function. The vector $\mathbf y$ denotes the velocity field, $p$ is the pressure, $D\mathbf{\mathbf{y}}=
\frac{1}{2}\big (\nabla \mathbf{\mathbf{y}}+(\nabla \mathbf{\mathbf{y}})^T\big )$ is the symmetric part of the velocity gradient, $\mathbf{u}$ is the given body force
 and $\Omega$ is an open bounded subset of $\RR^n$  ($n=2$ or $n=3$) .\par

System  (\ref{equation_etat}) can be used to model  steady incompressible electro-rheological fluids. It is based on the assumption that electro-rheological materials, composed by suspensions of particles in a fluid, can be considered as a homogenized single continuum media. The corresponding viscosity has the property of switching between shear-thinning and shear-thickening behavior under the application of a magnetic field. This model is described  and analyzed in \cite{Raj}, \cite{Ruz}, \cite{Ruz1} or \cite{Din}. More recently, in  \cite{Crispo1}, the authors proved the existence and uniqueness of a $C^{1,\gamma}(\bar{\Omega})\cap\mathbf W^{2,2}({\Omega})$ solution under smallness data conditions for system \ref{equation_etat}. This regularity result was the motivation for the analysis of the associated distributed optimal control problem that we describe in what follows.

Let us  look for the control $\mathbf{u}$ and the corresponding  $\mathbf{y}_u$  solution of $(\ref{equation_etat})$ such that the pair $(\mathbf{u},\mathbf{y}_u)$ solves  
%************************************************************%************************************************************
\vspace{0.3cm}
\begin{equation}\label{controlproblem}
(P_{\alpha})\left\{
  \begin{array}{ll}
    \mbox{Minimize}\,\, J(\mathbf{u},\mathbf{y}_u)\vspace{0.3cm}\\
    \mbox{subject to}\,\, (\ref{equation_etat})
  \end{array}
\right.\vspace{0.3cm}
\end{equation}
%************************************************************%************************************************************
\noindent where  the functional $J:\mathbf L^2(\Omega)\times\mathbf W^{1,2}_0(\Omega)\rightarrow\RR$ is given by
%************************************************************%************************************************************
\vspace{0.3cm}
\begin{equation}\label{costfuncional}
J(\mathbf{u},\mathbf{y})=\frac{1}{2}\int\limits_{\Omega}|\mathbf{\mathbf{y}}_u-\mathbf{\mathbf{y}}_d|^2\,dx+\frac{\nu}{2}\int\limits_{\Omega}|\mathbf{u}|^2\,dx.\vspace{0.3cm}
\end{equation}
%************************************************************%************************************************************
\noindent and $\mathbf{\mathbf{y}}_d$ denotes a fixed element of $\mathbf L^2(\Omega)$.\par
  
Such type of optimal control problems has been a subject of intensive research in the past decades. For  non-Newtonian fluid equations we mention the results in  \cite{Arada}, \cite{Casas1}, \cite{Casas2}, \cite{TG}, \cite{Gunz}, \cite{Slawig} and \cite{Wach} where the authors used several techniques to deal properly with the shear-thinning and shear-thickening viscosity laws, defined both in 2D and 3D domains. For the existence of solution, such techniques consist in exploring correctly the properties of the tensor $S$ in order to establish compactness results necessary for the application of the direct method of the Calculus of Variations. 
Our purpose here is to show that, based on the regularity results obtained in \cite{Crispo1}, it is possible to  easily extend those techniques to the case of electro-rheological fluids modeled by (\ref{equation_etat}). Treating the optimality conditions associated to problem $(P_{\alpha})$ is also an important, yet delicate, issue. We will therefore study this problem in a forthcoming work. 

In section 2 we introduce the notation that we are going to use, and  recall some useful results. In section 3 we characterize  the tensor $S$ including its continuity, coercivity and monotonicity properties. Finally, in section 4 we prove the main existence result of this paper.

\section{Notation and classical results}

We denote by ${\cal D}(\Omega)$ the space of infinitely differentiable functions with compact support in $\Omega$, ${\cal D}'(\Omega)$ denotes its dual (the space of distributions).
The standard Sobolev spaces are represented by $\mathbf W^{k,\alpha}(\Omega)$ ($k\in \NN$ and
$1<\alpha<\infty$), and their norms by $\|\cdot\|_{k,p}$. We set
$\mathbf W^{0,\alpha}(\Omega)\equiv \mathbf L^{\alpha}(\Omega)$ and $\|\cdot\|_{\alpha}\equiv
\|\cdot\|_{L^{\alpha}}$. The dual space of  $\mathbf W^{1,\alpha}_0(\Omega)$ is denoted by
 $\mathbf W^{-1,\alpha '}(\Omega)$ and its norm by $\|\cdot\|_{-1,\alpha '}$. We consider the space of divergence free functions defined by

        $$ {\cal V}=\big\{\psi\in {\cal D}(\Omega)
        \mid \nabla\cdot \psi=0\big\},  $$

\noindent to eliminate the pressure in the weak formulation. The space $\mathbf V_{\alpha}$ is the closure of ${\cal V}$ with respect to the gradient norm, i.e.

	$${\mathbf V}_{\alpha}=
        \left\{\psi\in \mathbf W_0^{1,\alpha}(\Omega)
        \mid \nabla\cdot \psi=0\right\}. $$

\noindent The space of H\"{o}lder continuous functions is a Banach space defined as

$$C^{m,\gamma}(\bar{\Omega})\equiv\{\mathbf y\in C^m(\bar{\Omega}):\|\mathbf y\|_{C^{m,\gamma}(\bar{\Omega})}<\infty\}  $$

\noindent where

\begin{equation}\label{normC}
\|\mathbf{\mathbf{y}}\|_{C^{m,\gamma}(\bar{\Omega})}\equiv \sum_{|\alpha|=0}^m\|D^{\alpha}\mathbf{\mathbf{y}}\|_{\infty}+[\mathbf{\mathbf{y}}]_{C^{m,\gamma}(\bar{\Omega})},        \vspace{0.3cm}
\end{equation}

\noindent and the semi-norm is denoted by

$$[\mathbf{y}]_{C^{m,\gamma}}(\bar{\Omega})\equiv\sum_{|\alpha|=m}\sup_{\{x_1,x_2\in\bar{\Omega}, x_1\neq x_2\}}\frac{|D^{\alpha}\mathbf{\mathbf{y}}(x_1)-D^{\alpha}\mathbf{\mathbf{y}}(x_2)|}{|x_1-x_2|^{\gamma}}<+\infty,  $$

\noindent for $m$ a nonnegative integer,  $0<\gamma<1$ and where

$$D^{\alpha}\mathbf y\equiv\frac{\partial^{|\alpha|}\mathbf y}{\partial x_1^{\alpha_1}...\partial x_n^{\alpha_n}}        \vspace{0.3cm},$$

\noindent with $\alpha=(\alpha_1,...,\alpha_n)$, $\alpha_i\in\NN_0$ and $|\alpha|=\sum\limits_{i=1}^n\alpha_i$.\\

Now we recall  two fundamental classical inequalities.

\begin{lemma}[Poincar\'{e}'s inequality]
Let $\mathbf y\in\mathbf W^{1,\alpha}_0(\Omega)$ with $1\leq\alpha<+\infty$. There exists a constant $C_1$ depending on $\alpha$ and $\Omega$ such that
\vspace{0.3cm}
$$\|\mathbf{\mathbf{y}}\|_{\alpha}\leq C_1(\alpha,\Omega)\|\nabla \mathbf{\mathbf{y}}\|_{\alpha}.\vspace{0.3cm}$$
\end{lemma}

\noindent {\it Proof.}
Found in \cite{Brezis}.

\begin{lemma}[Korn's inequality]
Let $\mathbf y\in\mathbf W^{1,\alpha}_0(\Omega)$ with $1<\alpha<+\infty$. There exists a constant $C_2$ depending on $\Omega$ such that
\vspace{0.3cm}
$$C_2(\Omega)\|\mathbf{\mathbf{y}}\|_{1,\alpha}\leq\|D\mathbf y\|_{\alpha}.\vspace{0.3cm}$$
\end{lemma}

\noindent {\it Proof.}  Found in \cite{Malek1}.

\vspace{0.3cm}

Finally two simple, yet very useful, properties of the convective term. 

\begin{lemma}
Let us consider $\mathbf u$ in $\mathbf{V}_{2}\,$ and $\,\mathbf{v}$, $\mathbf w$ in $\mathbf W^{1,2}_{0}(\Omega)$. Then
\vspace{0.3cm}
\begin{equation}\label{sym}
(\mathbf u \cdot\nabla \mathbf v,\mathbf w)=-(\mathbf u \cdot\nabla \mathbf w,\mathbf v)\quad\mbox{and}\quad (\mathbf u\cdot\nabla \mathbf v, \mathbf v)=0.
\end{equation}
\end{lemma}

\section{Properties of the extra tensor $S$}

Denoting by $\RR^{n\times n}_{sym}$ the set of all symmetric $n\times n$ matrices, we assume that  tensor $S: \RR^{n\times n}_{sym}\longrightarrow \RR^{n\times n}_{sym}$
has a potential, i.e.,  there exists a function $\Phi\in C^2(\RR^+_0,\RR^+_0)$ with $\Phi(0)=0$ such that

	\begin{eqnarray*}
	S_{ij}(\eta)=
	\frac{\partial \Phi(|\eta|^2)}{\partial\eta_{ij}}=
	2\Phi'(|\eta|^2)\,\eta_{ij}, \qquad S(0)=0 
	\end{eqnarray*}

\noindent for all $\eta\in \RR^{n\times n}_{sym}$. An example of such tensor is the one we are going to work with, namely

$$S(\eta)=(1+|\eta|)^{\alpha(x)-2}\eta,$$

\noindent where $\alpha(x)$ is a continuous function in $\bar{\Omega}$ such that

$$\alpha(x):\bar{\Omega}\rightarrow(1,+\infty)$$

\noindent and

\begin{equation}\label{a0}
1<\alpha_0\leq\alpha(x)\leq\alpha_{\infty}<+\infty
\end{equation}

$$\min\,\alpha(x)=\alpha_0,$$

$$\max\,\alpha(x)\leq\alpha_{\infty}, \quad\mbox{for}\quad\alpha_{\infty}>2.$$

The case $1< \alpha(x) <2$ corresponds to shear-thinning viscosity fluids, while $\alpha(x) >2$ corresponds to shear-thickening fluids. The case $\alpha(x)=2$ corresponds to a Newtonian fluid.
For such function $\alpha(x)$, it can be proved that tensor $S$ satisfies standard properties.

\begin{proposition}
Consider $\alpha(x) \in(1,\infty)$, $C_3$ and $C_4$ positive constants. Then the following inequalities hold

\begin{description}
 \item {\bf A1\,-} $$\left|\frac{\partial S_{k\ell}(\eta)}{\partial \eta_{ij}}\right|\leq C_3
	\left(1+|\eta|\right)^{\alpha(x)-2}.$$

 \item {\bf A2\,-} $$S '(\eta):\zeta:\zeta=\sum_{ijk\ell}
	\textstyle\frac{\partial S_{k\ell}(\eta)}{\partial \eta_{ij}}
	\zeta_{k\ell}\zeta_{ij}\geq C_4
	 \left(1+|\eta|\right)^{\alpha(x)-2}|\zeta|^2$$
\noindent for all $\eta, \zeta\in \RR^{n\times n}_{sym}$ and $i,j,k,\ell=1,\cdots, d$.
\end{description}
\end{proposition}

\noindent {\it Proof.}
In fact,

\begin{eqnarray}\label{a1}
\nonumber\left|\frac{\partial S_{kl}}{\partial\eta_{ij}}\right|&=&\left|(\alpha(x)-2)\left(1+|\eta|\right)^{\alpha(x)-3}\frac{\eta_{ij}}{|\eta|}\eta_{kl}+\left(1+|\eta|\right)^{\alpha(x)-2}\delta_{ik}\delta_{jl}\right|\\ 
&\leq&\left|\alpha(x)-2\right|\left(1+|\eta|\right)^{\alpha(x)-3}\frac{|\eta_{ij}\eta_{kl}|}{|\eta|}+\left(1+|\eta|\right)^{\alpha(x)-2}\left|\delta_{ik}\delta_{jl}\right|.
\end{eqnarray}

\noindent Taking into account that

  $$\delta_{ik}\delta_{jl}=\left\{\begin{array}{cc}
    1 & \mbox{if}\,\, i=k,j=l \\ 
    0 & \mbox{otherwise} \\ 
      \end{array}\right.$$
      
 \noindent and 
 
 $$|\eta_{ij}\eta_{kl}|\leq|\eta|^2,$$     

\noindent we can write

\begin{eqnarray}\label{a2}
\nonumber(\ref{a1})&\leq&\left|\alpha(x)-2\right|\left(1+|\eta|\right)^{\alpha(x)-3}|\eta|+\left(1+|\eta|\right)^{\alpha(x)-2}\\
\nonumber&\leq&\left|\alpha(x)-2\right|\left(1+|\eta|\right)^{\alpha(x)-3}(1+|\eta|)+\left(1+|\eta|\right)^{\alpha(x)-2}\\
\nonumber&=&\left|\alpha(x)-2\right|\left(1+|\eta|\right)^{\alpha(x)-2}+\left(1+|\eta|\right)^{\alpha(x)-2}\\
&=&(\left|\alpha(x)-2\right|+1)\left(1+|\eta|\right)^{\alpha(x)-2}.
\end{eqnarray}

\noindent If $\alpha(x)-2\geq 0$, we have

\begin{eqnarray*}
(\ref{a2})&=&(\alpha(x)-1)\left(1+|\eta|\right)^{\alpha(x)-2}\\
&\leq&(\alpha_{\infty}-1)\left(1+|\eta|\right)^{\alpha(x)-2}.\\
\end{eqnarray*}

\noindent Otherwise, if $\alpha(x)-2< 0$, we have

\begin{eqnarray*}
(\ref{a2})&=&(3-\alpha(x))\left(1+|\eta|\right)^{\alpha(x)-2}\\
&\leq&(3-\alpha_0)\left(1+|\eta|\right)^{\alpha(x)-2},\\
\end{eqnarray*}

\noindent and therefore we have 

$$\left|\frac{\partial S_{k\ell}(\eta)}{\partial \eta_{ij}}\right|\leq\left\{  \begin{array}{ccc}
   (3-\alpha_{0})\left(1+|\eta|\right)^{\alpha(x)-2}& \mbox{if} &  \alpha(x)-2< 0\\ 
    (\alpha_{\infty}-1)\left(1+|\eta|\right)^{\alpha(x)-2} & \mbox{if}& \alpha(x)-2\geq 0 .\\ 
  \end{array}\right.$$

\noindent This proves inequality $\bf A1.$ 

\vspace {0.5 cm}

\noindent In order to obtain $\bf A2$, we write

\begin{eqnarray}\label{a3}
\nonumber && S '(\eta):\zeta:\zeta = 
\sum_{ijk\ell}
	\textstyle\frac{\partial S_{k\ell}(\eta)}{\partial \eta_{ij}}
	\zeta_{k\ell}\zeta_{ij}\\
\nonumber&=&\sum_{ijk\ell}\left[(\alpha(x)-2)\left(1+|\eta|\right)^{\alpha(x)-3}\frac{\eta_{ij}\eta_{kl}}{|\eta|} +\left(1+|\eta|\right)^{\alpha(x)-2}\delta_{ik}\delta_{jl}\right]\zeta_{ij}\zeta_{k\ell}\\
&=&(\alpha(x)-2)\frac{\left(1+|\eta|\right)^{\alpha(x)-3}}{|\eta|}\sum_{ijk\ell}\eta_{ij}\eta_{kl}\zeta_{ij}\zeta_{k\ell}+\nonumber\left(1+|\eta|\right)^{\alpha(x)-2}\sum_{ijk\ell}\delta_{ik}\delta_{jl}\zeta_{ij}\zeta_{k\ell}.\\
\end{eqnarray}

\noindent Considering that

$$\sum_{ijk\ell}\eta_{ij}\eta_{kl}\zeta_{ij}\zeta_{k\ell}=\sum_{ij}\eta_{ij}\zeta_{ij}\sum_{kl}\eta_{kl}\zeta_{kl}=|\eta:\zeta|^2$$

\noindent and

$$\sum_{ijk\ell}\delta_{ik}\delta_{jl}\zeta_{ij}\zeta_{k\ell}=\sum_{ij}\zeta_{ij}\zeta_{ij}=|\zeta|^2,$$

\noindent  expression (\ref{a3}) is equal to

\begin{eqnarray}\label{a4}
&&(\alpha(x)-2)\frac{\left(1+|\eta|\right)^{\alpha(x)-3}}{|\eta|}|\eta:\zeta|^2+\left(1+|\eta|\right)^{\alpha(x)-2}|\zeta|^2.
\end{eqnarray}

\noindent Taking into account that $\alpha(x)-2< 0$, $\alpha_{0}\leq\alpha(x)$ and $|\eta:\zeta|^2\leq|\eta|^2|\zeta|^2$, it follows 

\begin{eqnarray*}
(\ref{a4})&\geq&(\alpha(x)-2)\frac{\left(1+|\eta|\right)^{\alpha(x)-3}}{|\eta|}|\eta|^2|\zeta|^2+\left(1+|\eta|\right)^{\alpha(x)-2}|\zeta|^2\\
&=&\left((\alpha(x)-2)(1+|\eta|)^{\alpha(x)-3}|\eta|+(1+|\eta|)^{\alpha(x)-2}\right)|\zeta|^2\\
&\geq&\left((\alpha(x)-2)(1+|\eta|)^{\alpha(x)-3}(1+|\eta|)+(1+|\eta|)^{\alpha(x)-2}\right)|\zeta|^2\\
&=&\left((\alpha(x)-2)(1+|\eta|)^{\alpha(x)-2}+(1+|\eta|)^{\alpha(x)-2}\right)|\zeta|^2\\
&=&(\alpha(x)-1)(1+|\eta|)^{\alpha(x)-2}|\zeta|^2\\
&\geq&(\alpha_{0}-1)(1+|\eta|)^{\alpha(x)-2}|\zeta|^2
\end{eqnarray*}

\noindent Instead, $\alpha(x)-2\geq 0$ gives

\begin{eqnarray*}
(\ref{a4})&\geq&\left(1+|\eta|\right)^{\alpha(x)-3}\left((\alpha(x)-2)\frac{|\eta:\zeta|^2}{|\eta|}+(1+|\eta|)|\zeta|^2\right)\\
&\geq&\left(1+|\eta|\right)^{\alpha(x)-3}(1+|\eta|)|\zeta|^2\\
&=&\left(1+|\eta|\right)^{\alpha(x)-2}|\zeta|^2.
\end{eqnarray*}

\noindent Then we have

$$S '(\eta):\zeta:\zeta\geq\left\{  \begin{array}{ccc}
    (\alpha_{0}-1)\left(1+|\eta|\right)^{\alpha(x)-2}|\zeta|^2& \mbox{if} &  \alpha(x)-2< 0\\ 
    \left(1+|\eta|\right)^{\alpha(x)-2}|\zeta|^2 & \mbox{if}& \alpha(x)-2\geq 0 .\\ 
  \end{array}\right.$$
  
  \noindent This proves inequality $\bf A2.$ 

\vspace {0.5 cm}

  \noindent Assumptions ${\bf A1}$-${\bf A2}$ imply the following standard properties for $S$ (see \cite{Malek1}).\\

\begin{proposition} In the same conditions we have

\begin{description}
\item [1.] Continuity
\begin{equation}\label{conti}
|S(\eta)|\leq (1+|\eta|)^{\alpha(x)-2}|\eta|,
\end{equation} 

\item [2.] Coercivity
\begin{equation}\label{coerc}
S(\eta):\eta\geq \left\{ \begin{array}{ccc}

   \nu(1+|\eta|)^{\alpha(x) -2}|\eta|^2, & \mbox{if} & \alpha(x)-2<0 \\ 
    |\eta|^2, & \mbox{if} & \alpha(x)-2\geq0  \\ 
  \end{array}\right.
\end{equation}

\item [3.] Monotonicity
\begin{equation}\label{mono}
(S(\eta)-S(\zeta)):(\eta-\zeta)\geq\nu (1+|\eta|+|\zeta|)^{\alpha(x) -2}|\eta-\zeta|^2.
\end{equation}
\end{description}
\end{proposition}

\noindent {\it Proof.}
Continuity is trivially derived from

\begin{eqnarray*}
|S(\eta)|= \left|(1+|\eta|)^{\alpha(x)-2}\eta\right |=\left|(1+|\eta|)^{\alpha(x)-2}\right | |\eta|=(1+|\eta|)^{\alpha(x)-2} |\eta|.
\end{eqnarray*} 

\noindent Coercivity is equivalent to monotonocity taking $S(\zeta)=\zeta=0_M.$ Therefore it is enough to prove monotonocity. Taking into account that

\begin{eqnarray*}
S_{ij}(\eta)-S_{ij}(\zeta)&=&\int\limits\limits_{0}^1\frac{\partial}{\partial t} S_{ij}(t\eta+(1-t)\zeta)\,dt\\
&=&\int\limits\limits_{0}^1\sum\limits_{kl}\,\frac{\partial S_{ij}(t\eta+(1-t)\zeta)}{\partial D_{kl}}(\eta-\zeta)_{kl} 
\end{eqnarray*}

\noindent we can write

\begin{eqnarray}\label{a5}
\nonumber (S(\eta)-S(\zeta):(\eta-\zeta)&=&\int\limits\limits_{0}^1\sum\limits_{ij}\,\sum\limits_{kl}\frac{\partial S_{ij}(t\eta+(1-t)\zeta)}{\partial D_{kl}}(\eta-\zeta)_{kl} :(\eta-\zeta)_{ij} \,dt\\
&=&\int\limits\limits_0^1S'(t\eta+(1-t)\zeta)):(\eta-\zeta):(\eta-\zeta)\,dt
\end{eqnarray}

\noindent Using $\bf{A}2$ and considering $\alpha(x)-2< 0$ we have

\begin{eqnarray}\label{ai6}
(\ref{a5})\geq\int\limits\limits_{0}^1(\alpha_{0}-1)\left(1+|t\eta+(1-t)\zeta|\right)^{\alpha(x)-2}|\eta-\zeta|^2\,dt
\end{eqnarray}

\noindent Since $t\in[0,1]$ then

\begin{eqnarray}
1+|t\eta+(1-t)\zeta|\leq 1+|\eta+\zeta|\leq 1+|\eta|+|\zeta|
\end{eqnarray}

\noindent and therefore, we can write

\begin{eqnarray*}
(\ref{ai6})\geq\int\limits\limits_{0}^1(\alpha_{0}-1)\left((1+|\eta|+|\zeta|)\right)^{\alpha(x)-2}|\eta-\zeta|^2\,dt
\end{eqnarray*}

\noindent Using $\bf{A}2$ and considering $\alpha(x)-2\geq 0$ we show that

\begin{eqnarray}\label{a6}
\nonumber(\ref{a5})&\geq&\int\limits\limits_{0}^1\left(1+|t\eta+(1-t)\zeta|\right)^{\alpha(x)-2}|\eta-\zeta|^2\,dt\\
\nonumber&\geq&\int\limits\limits_{0}^1 1^{\alpha(x)-2}|\eta-\zeta|^2\,dt\\
&\geq& |\eta-\zeta|^2.
\end{eqnarray}

\noindent Then we conclude 

 $$(S(\eta)-S(\zeta)):(\eta-\zeta)\geq \left\{ \begin{array}{ccc}
   (\alpha_{0}-1)(1+|\eta|+|\zeta|)^{\alpha(x) -2}|\eta-\zeta|^2 & \mbox{if} & \alpha(x)-2<0 \\ 
    |\eta-\zeta|^2 & \mbox{if} & \alpha(x)-2\geq0.  \\ 
  \end{array}\right.$$

\section{Main Result}

\begin{definition}\label{def2}
Assume that $\mathbf{u}\in\mathbf L^2(\Omega)$. A function $\mathbf{y}$ is a $C^{1,\gamma}$-solution of (\ref{equation_etat}) if $\mathbf{y}\in C^{1,{\gamma}}(\bar{\Omega})$, for $\gamma\in (0,1)$, $div\, \mathbf{y}=0$, $\mathbf{y}|_{\partial\Omega}=0$ and it satisfies the following  equality
\vspace{0.3cm}
\begin{equation}\label{strongsol}
(S(D\mathbf{y}),D\varphi)+(\mathbf{y}\cdot\nabla \mathbf{y},\varphi)=(\mathbf{u},
	\varphi),\quad \mbox{for all}\,\,  \varphi\in \mathbf{V}_{2},
	\end{equation}
	
\noindent  where $(\cdot,\cdot)$ denotes the scalar product in $L^2	(\Omega).$
\end{definition}

Next proposition, due to \cite{Crispo1},  presents an existence and uniqueness result of a $C^{1,\gamma_0}$ solution of system (\ref{equation_etat}) with certain conditions imposed to $\mathbf{u}$, but with no aditional conditions on the exponent $\alpha$.

\begin{proposition}\label{theo1}
We assume that $\mathbf u \in \mathbf{L}^{q}(\Omega),$ for some $q>n$. Let us consider $\Omega$ a $C^{1,\gamma_0}$ domain, and $\alpha\in C^{0,\gamma_0}(\bar{\Omega})$, with $\gamma_0=1-\frac{n}{q}$. Then, for any $\gamma<\gamma_0$, there exist positive constants $C_5$ and $C_6$, depending on $\|\alpha\|_{C^{0,\gamma}(\bar{\Omega})}$, $n$, $q$ and $\Omega$ such that, if $\|\mathbf{u}\|_q< C_5$,   there exists a $C^{1,\gamma}$ solution $(\mathbf{y},p)$ of problem (\ref{equation_etat}) verifying

\begin{equation}\label{estimative1}
\|\mathbf{y}\|_{C^{1,\gamma}(\bar{\Omega})}+\|p\|_{C^{0,\gamma}(\bar{\Omega})}\leq C_6\|\mathbf{u}\|_q.
\end{equation}

\noindent Furthermore, there exists a constant $C_7$ depending on $\alpha_{0}$, $\|\alpha\|_{C^{0,\gamma}(\bar{\Omega})}$, $n$, $q$ and $\Omega$ such that, if $\|\mathbf{u}\|_q\leq C_7$, the solution is unique.
\end{proposition}

\begin{proposition}\label{theo2}
Assume that properties $\bf A1$ and $\bf A2$ are fulfilled. Considering ${\mathbf{y}\in C^{1,\gamma}(\bar{\Omega})}$ we have

\begin{equation}\label{t2}
\|D\mathbf{y}\|_2\leq C_8\|\mathbf{u}\|_2,
\end{equation}

\noindent where $\mathbf{y}$ is the associated state to $\mathbf{u}$.
\end{proposition}

\noindent {\it Proof.} Taking $\varphi=\mathbf{y}$ in (\ref{strongsol}) and recalling the convective term properties, we have

\begin{equation}\label{t1}
(S(D\mathbf{y}),D\mathbf{y})=(\mathbf{u},
	\mathbf{y}).
	\end{equation}

\noindent Since ${\mathbf{y}\in C^{1,\gamma}(\bar{\Omega})}$ then $\mathbf{y}$  belongs to $C^1(\bar{\Omega})$ which means that $\mathbf{y}$ and $D\mathbf{y}$ are bounded functions in $\bar{\Omega}$ and consequently belong to $\mathbf{L}^{\alpha}(\Omega),$ for any $\alpha >1$.  In particular we consider $\mathbf y\in\mathbf L^2(\Omega)$. Then, by using H\"{o}lder's inequality and the Poincar\'{e} and Korn inequalities there exists a constant $C_8$ such that

\begin{eqnarray*}
|(\mathbf{u},
	\mathbf{y})|\leq\|\mathbf{u}\|_2\|\mathbf{y}\|_2\leq C_8 \|\mathbf{u}\|_2\|D\mathbf{y}\|_2.
	\end{eqnarray*}

\noindent  On the other hand, by coercivity  we write

$$\|D\mathbf{y}\|^2_2\leq (S(D\mathbf{y}),D\mathbf{y}).$$

\noindent Putting together both inequalities with (\ref{t1}) we prove the pretended result.\par

\vspace {0.3 cm}

Once we have the guarantee of existence of a solution of problem (\ref{equation_etat}) provided by Proposition \ref{theo1} and an estimative of energy for $D\mathbf{y}$ given by proposition \ref{theo2}, we can now formulate and prove the following existence result for the control problem $(P_{\alpha})$.

\begin{theorem}[Main Result]\label{MainResult}
Assume that $\bf A1$-$\bf A2$ are fulfilled, with ${1<\alpha\leq 2}$. Then $(P_{\alpha})$ admits at least one solution.
\end{theorem}

To prove this theorem we need to establish some important results.

\begin{proposition}\label{prop7}
Assume that $(\mathbf{u}_k)_{k>0}$ converges to $\mathbf{u}$ weakly in $\mathbf L^2(\Omega)$. Let $\mathbf y_k$ be the associated state to $\mathbf u_k$. Then there exists  $\mathbf{y}\in\mathbf{W}^{1,2}_0(\Omega)$ and $\tilde{S}\in\mathbf{L}^2(\Omega)$ such that the following convergences are verified

\begin{equation}\label{t5}
(\mathbf{y}_k)_k\rightharpoonup \mathbf y\quad\mbox{in}\quad\mathbf{W}^{1,2}_0(\Omega)
\end{equation}

\begin{equation}\label{t4}
(D\mathbf{y}_k)_k\rightharpoonup D\mathbf y\quad\mbox{in}\quad\mathbf{L}^2(\Omega)
\end{equation}

\begin{equation}\label{t6}
(S(D\mathbf{y}_k))_k\rightharpoonup \tilde{S}\quad\mbox{in}\quad\mathbf{L}^2(\Omega).
\end{equation}

\end{proposition}

\noindent {\it Proof.}
The convergence of $(\mathbf{u}_k)_{k>0}$ to $\mathbf{u}$ in the weak topology of $\mathbf L^2(\Omega)$ implies that $(\mathbf{u}_k)_{k>0}$ is bounded, i.e, there exists a positive constant $M$ such that

\begin{equation}\label{t3}
\|\mathbf{u}_k\|_2\leq M,\quad\mbox{for}\,\,k>k_0.\vspace{0.3cm}
\end{equation}

\noindent Due to (\ref{t2}) and (\ref{t3}), it follows that

$$\|D\mathbf{y}_k\|_{2}\leq C_8 M.\vspace{0.3cm}$$

\noindent By  Korn's inequality $\mathbf{y}_k$ is then bounded in $\mathbf W^{1,2}_0(\Omega)$ and thus there is a subsequence still indexed in $k$ that weakly converges to a certain $\mathbf{y}$ in $\mathbf W^{1,2}_0(\Omega)$. Moreover, by using a Sobolev's compact injection, $\mathbf{y}_k$ converges strongly (then also weakly) to $\mathbf{y}$ in $\mathbf{L}^2(\Omega)$. It is straightforward to conclude (\ref{t4}).

\noindent Finally, the previous estimate, together with $(\ref{conti})$ implies

\begin{eqnarray*}{\label{desig1}}
\|S(D\mathbf{y}_k)\|^{2}_{2}&\leq& \int\limits_{\Omega}(1+|D\mathbf{y}_k|)^{(\alpha(x) -2)2}|D\mathbf{y}_k|^{2}\,dx\\
&\leq& \int\limits_{\Omega}(1+|D\mathbf{y}_k|)^{(\alpha(x) -2)2}(1+|D\mathbf{y}_k|)^2\,dx\\
&= &\int\limits_{\Omega}(1+|D\mathbf{y}_k|)^{2(\alpha(x) -1)}\,dx\\
&\leq& C_8\int\limits_{\Omega}(1+|D\mathbf{y}_k|^{2(\alpha(x) -1)})\,dx\\
&\leq& C_8 \left(|\Omega|+\int\limits_{\Omega}|D\mathbf{y}_k|^{2(\alpha_{\infty} -1)}\,dx\right)\\
&=& C_8 \left(|\Omega|+\|D\mathbf{y}_k\|^{2(\alpha_{\infty}  -1)}_{2(\alpha_{\infty} -1)}\right).
\end{eqnarray*}

\noindent This last expression is  bounded once $D\mathbf{y}_k\in C(\bar{\Omega})$ and consequently the sequence $(S(D\mathbf y_k))_k$ is bounded in $\mathbf L^{2}(\Omega)$ and we finish the proof by establishing the existence of a subsequence, still indexed by $k$, and $\tilde{S}\in\mathbf L^{2}(\Omega)$ such that  $(S(D\mathbf{y}_k))_{k>0}$ weakly converges to $\tilde{S}\in\mathbf L^{2}(\Omega)$. 

\begin{proposition}\label{prop7}
Assume that (\ref{t4}), (\ref{t5}) and (\ref{t6}) are verified. Then the weak limit of $(\mathbf{y}_k)_k$, $\mathbf{y}$, is the solution of (\ref{strongsol}) corresponding to $\mathbf{u}\in\mathbf L^2(\Omega)$.\end{proposition}

\noindent {\it Proof.}
Let us consider

\begin{equation}\label{desig4}
(S(D\mathbf{y}_k)-S(\mathbf y),D\varphi)+(\mathbf{y}_k\cdot\nabla \mathbf{y}_k-\mathbf y\cdot\nabla \mathbf y,\varphi)=(\mathbf{u}_k-\mathbf{u},\varphi),
\end{equation}

\noindent for all $\varphi\in\textbf{V}_{2}.$ Taking into account the convective term properties and the regularity results assumed on $\mathbf{y}$, we have

\begin{eqnarray*}
|(\mathbf{y}_k\cdot\nabla \mathbf{y}_k-\mathbf y\cdot\nabla \mathbf y,\varphi)
=|((\mathbf{y}_k-\mathbf y)\cdot\nabla \mathbf{y}_k,\varphi)+(\mathbf y\cdot\nabla(\mathbf{y}_k-\mathbf y),\varphi)|\\
 =|((\mathbf y_k-\mathbf y)\cdot\nabla\mathbf  y_k,\varphi)-(\mathbf y\cdot\nabla \varphi,(\mathbf y_k-\mathbf y))|\\
\leq|((\mathbf y_k-\mathbf y)\cdot\nabla \mathbf y_k,\varphi)|+|(\mathbf y\cdot\nabla \varphi,(\mathbf y_k-\mathbf y))|\\
\leq C_E^2\left(\|\nabla \mathbf y_k\|_{2}\|\varphi\|_{4}+\|\mathbf y\|_{4}
\|\nabla\varphi\|_{2}\right)\|\mathbf y_k-\mathbf y\|_{4}
\rightarrow 0, \quad\mbox{when}\quad k\rightarrow + \infty.
\end{eqnarray*}

\noindent This result is a consequence of the compact injection of $\mathbf W^{1,2}_0(\Omega)$ into $\mathbf L^4(\Omega)$ which provide a strong convergence in $\mathbf L^4(\Omega)$ due to  (\ref{t5}). Note that $C_E$ corresponds to the embedding constant. 

\noindent Hence, passing to the limit in

$$(S(D\mathbf y_k),D\varphi)+(\mathbf y_k\cdot\nabla \mathbf y_k,\varphi)=(\mathbf u_k,\varphi),\quad\mbox{for all}\quad\varphi\in\mathbf V_{2},\vspace{0.3cm}$$

\noindent we obtain

\begin{equation}\label{eq1}
(\tilde{S},D\varphi)+(\mathbf y\cdot\nabla \mathbf y,\varphi)=(\mathbf u,\varphi),\quad\mbox{for all}\quad\varphi\in\mathbf V_{2},\vspace{0.3cm}
\end{equation}

\noindent In particular, taking $\varphi=\mathbf y$ and considering $(\ref{sym})$ we may write

\begin{equation}\label{eq2}
(\tilde{S},D\mathbf y)=(\tilde{S},D\mathbf y)+(\mathbf y\cdot\nabla \mathbf y,\mathbf y)=(\mathbf u,\mathbf y).
\end{equation}

\noindent On the other hand, the monotonocity assumption $(\ref{mono})$ gives

\begin{equation}\label{eq3}
(S(D\mathbf y_k)-S(D\varphi),D(\mathbf y_k)-D\varphi)\geq 0,\quad\mbox{for all}\quad\varphi\in\mathbf V_{2}.\vspace{0.3cm}
\end{equation}

\noindent Since,

$$(S(D\mathbf y_k),D\mathbf y_k)=(\mathbf u_k,\mathbf y_k),\vspace{0.3cm}$$
 replacing the first member in $(\ref{eq3})$, we obtain

$$(\mathbf u_k,\mathbf y_k)-(S(D \mathbf y_k),D\varphi)-(S(D\varphi), D\mathbf y_k -D\varphi)\geq 0, \quad \mbox{for all}\quad\varphi\in\mathbf V_{2}.$$

\noindent Passing to the limit it follows

$$(\mathbf u,\mathbf y)-(\tilde{S},D\varphi)-(\tau(D\varphi), D\mathbf y-D\varphi)\geq 0, \quad\mbox{for all}\quad\varphi\in\mathbf V_{2}.$$

\noindent This inequality together with $(\ref{eq2})$, implies that

$$(\tilde{S}-S(D\varphi),D\mathbf y-D\varphi)\geq 0, \quad\mbox{for all}\quad\varphi\in\mathbf V_{2}$$

\noindent Taking $\varphi=\mathbf y-\lambda \mathbf v$ (see \cite{Lions1}), which is possible considering any $\mathbf{v}\in\mathbf V_2$ and $\lambda>0$, we have

\begin{equation}
(\tilde{S}-S(D(\mathbf y-\lambda \mathbf v)),D\mathbf y-D(\mathbf y-\lambda \mathbf v))\geq 0,\quad\mbox{for all}\quad\mathbf v\in\mathbf V_{2}
\end{equation}

\noindent which is equivalent to 

\begin{equation}
\lambda(\tilde{S}-S(D(\mathbf y-\lambda \mathbf v)),D\mathbf v))\geq 0,\quad\mbox{for all}\quad \mathbf v\in\mathbf V_{2}
\end{equation}

\noindent and then, since  $\lambda>0$, it comes

\begin{equation}
(\tilde{S}-S(D(\mathbf y-\lambda \mathbf v)),D\mathbf v))\geq 0,\quad\mbox{for all}\quad \mathbf v\in\mathbf V_{2}
\end{equation}

Passing to the limit when $\lambda\rightarrow 0$ and considering the continuity of $S$ we obtain

\begin{equation}
(\tilde{S}-S(D(\mathbf y)),D\mathbf v))\geq 0,\quad\mbox{for all}\quad \mathbf v\in\mathbf V_{2}.
\end{equation}

\noindent This implies that

$$\tilde{S}=S(D(\mathbf y))$$
and then

$$(S(D\mathbf y),D\varphi)+(\mathbf y\cdot\nabla \mathbf y,\varphi)=(\mathbf u,\varphi), \quad\mbox{for all}\quad\varphi\in\mathbf V_{2}.$$

\noindent Hence, $\mathbf{y}\equiv \mathbf y_u$, i.e, $\mathbf{y}$ is the solution associated to $\mathbf{u}$.

\begin{proposition}\label{prop7}
Assume that $\mathbf A_1$ and $\mathbf A_2$ are satisfied. Then $(\mathbf{y}_k)_k$ strongly converges to $\mathbf{y}_u$ in $\mathbf{W}^{1,2}_0(\Omega)$.
\end{proposition}

\noindent {\it Proof.}
Setting $\varphi=y_k-y_u$ in (\ref{strongsol}) and taking  (\ref{mono}) we obtain

\begin{equation}
( S(D\mathbf y_k)-S(D\mathbf y_u),D(\mathbf y_k-\mathbf y_u))\geq \|D(\mathbf y_k-\mathbf y_u)\|^2_{2}
\end{equation}

\noindent Therefore, using (\ref{sym}) and classical embedding results, we obtain

%\begin{eqnarray*}
%&&\limsup_{k\rightarrow +\infty} \|D(\mathbf y_k-\mathbf y_u)\|^2_{2}\\
%&\leq&\limsup_{k\rightarrow +\infty} (S(D\mathbf y_k)-S(D\mathbf y_u),D(\mathbf y_k-\mathbf y_u))\\
%&=&\limsup_{k\rightarrow +\infty}((\mathbf u_k-\mathbf u,\mathbf y_k-\mathbf y_u)-(\mathbf y_k\cdot\nabla \mathbf y_k-\mathbf y_u\cdot\nabla \mathbf y_u,\mathbf y_k-\mathbf y_u))\\
%&=&\limsup_{k\rightarrow +\infty}((\mathbf u_k-\mathbf u,\mathbf y_k-\mathbf y_u)-((\mathbf y_k-\mathbf y_u)\cdot\nabla \mathbf y_u,\mathbf y_k-\mathbf y_u))\\
%&\leq&\limsup_{k\rightarrow +\infty}\left((\mathbf u_k-\mathbf u,\mathbf y_k-\mathbf y_u)-\|\mathbf y_k-\mathbf y_u\|^2_{4}\|\nabla \mathbf y_u\|_{2}\right)\\
%&\leq&\limsup_{k\rightarrow +\infty}\left((\mathbf u_k-\mathbf u,\mathbf y_k-\mathbf y_u)-\|\mathbf y_k-\mathbf y_u\|^2_{4}\|\nabla \mathbf y_u\|_{2}\right)=0.
%\end{eqnarray*}

\begin{eqnarray*}
\|D(\mathbf y_k-\mathbf y_u)\|^2_{2}&\leq&(S(D\mathbf y_k)-S(D\mathbf y_u),D(\mathbf y_k-\mathbf y_u))\\
&=&((\mathbf u_k-\mathbf u,\mathbf y_k-\mathbf y_u)-(\mathbf y_k\cdot\nabla \mathbf y_k-\mathbf y_u\cdot\nabla \mathbf y_u,\mathbf y_k-\mathbf y_u))\\
&=&((\mathbf u_k-\mathbf u,\mathbf y_k-\mathbf y_u)-((\mathbf y_k-\mathbf y_u)\cdot\nabla \mathbf y_u,\mathbf y_k-\mathbf y_u))\\
&\leq&\left((\mathbf u_k-\mathbf u,\mathbf y_k-\mathbf y_u)-\|\mathbf y_k-\mathbf y_u\|^2_{4}\|\nabla \mathbf y_u\|_{2}\right)\rightarrow 0.
\end{eqnarray*}

\noindent and therefore, by Korn's inequality, $$\|\mathbf y_k- \mathbf y\|_{1,2}\rightarrow 0.$$

\vspace{0.3cm}

Now we can prove our main result.

\vspace{0.3cm}

\noindent {\it Proof of Theorem \ref{MainResult}.}
Let $(\mathbf u_k)_k$ be a minimizing sequence in $\mathbf{L}^2(\Omega)$ and $(\mathbf y_k)_k$ the sequence of associated states. Considering the properties of the functional $J$ defined by (\ref{costfuncional}), we obtain

$$\frac{\nu}{2}\|\mathbf u_k\|^2_2\leq J(\mathbf u_k, \mathbf y_k)\leq J(0, \mathbf y_0),\quad\mbox{for}\quad k>k_0$$

\noindent implying that $(\mathbf u_k)_k$ is bounded in $\mathbf L^2(\Omega)$. From Proposition $\ref{prop7}$, we deduce that $(y_k)$ converges strongly to $y_u$. $ J$ is a sum of quadratic terms and is convex. On the other hand, if 

$$(\mathbf{v}_k,\mathbf z_k)\rightarrow(\mathbf v, \mathbf z)\quad\mbox{ in}\quad\mathbf L^2(\Omega)\times \mathbf W^{1,2}_0(\Omega)$$ 
this implies
$$J(\mathbf{v}_k,\mathbf z_k)\rightarrow J(\mathbf v, \mathbf z)\quad\mbox{ in}\quad\RR$$ 
and then the  functional $J$ is also a continuous function. In fact, once we have 

\begin{eqnarray*} 
&&\left|J(\mathbf{v}_k, \mathbf{z}_k)-J(\mathbf{v}, \mathbf{z})\right|=\left|\|\mathbf z _k-\mathbf y_d\|^{2}_2+\|\mathbf v_k\|_2^2-\|\mathbf z-\mathbf y_d\|^2_2-\|\mathbf v\|_2^2\right|\\
&\leq& \left|\left( \|(\mathbf z_k-\mathbf y_d)-(\mathbf z-\mathbf y_d)\|_2+\|(\mathbf z-\mathbf y_d)\|_2\right)^2+\left(\|\mathbf v_k
-\mathbf v\|_2+\|\mathbf v\|_2\right)^2-\|\mathbf z-\mathbf y_d\|^2_2-\|
\mathbf v\|_2^2\right|\\
&\leq &\left|\left( \|(\mathbf z_k-\mathbf z\|_2+\|(\mathbf z-\mathbf y_d)\|_2\right)^2+\left(\|\mathbf v_k-\mathbf v\|_2+\|\mathbf v\|_2\right)^2-\|\mathbf z-\mathbf y_d\|^2_2-\|
\mathbf v\|_2^2\right|.
\end{eqnarray*}
%************************************************************%************************************************************
Since $\mathbf z_k \rightarrow \mathbf z $ strongly also in $\mathbf L^2(\Omega)$, the last expression converges to zero when $k\rightarrow\infty$ and therefore, $J$ is a semicontinuous function (see  \cite {Brezis}). We may now apply the direct method of the Calculus of Variations (see e.g. \cite{Dac})
%************************************************************%************************************************************

$$\inf_{k} J\leq J(\mathbf u,\mathbf y_u)\leq\liminf_{k}J(\mathbf u_k,\mathbf y_k)\leq\inf_{k} J$$
%************************************************************%************************************************************
to conclude that $(\mathbf u,\mathbf y_u)$ is in fact a minimizer and therefore a solution of the control problem $(P_{\alpha}).$

\begin{acknowledgements}
This work has been partially supported  by FCT (Portugal) through the
Research Centers CMA/FCT/UNL, CEMAT-IST,  grant SFRH/BPD/66638/2009 and the projects PTDC/MAT109973/2009,  EXCL/MAT-NAN/0114/2012.
\end{acknowledgements}

% BibTeX users please use one of
%\bibliographystyle{spbasic}      % basic style, author-year citations
%\bibliographystyle{spmpsci}      % mathematics and physical sciences
%\bibliographystyle{spphys}       % APS-like style for physics
%\bibliography{}   % name your BibTeX data base

% Non-BibTeX users please use

\end{document}